\documentclass[12pt]{amsart}
\usepackage{latexsym, amsmath, amssymb}
\usepackage{epsfig}
\usepackage{amscd}
\usepackage{pstricks,pst-node,cite}

\newtheorem{thm}{Theorem}[section]

\newtheorem{cor}[thm]{Corollary}
\newtheorem{exa}[thm]{Example}
\theoremstyle{remark}

\def\aa{{\mathcal A}}

\def\tt{{\mathcal T}}
\def\ttt{{\mathfrak t}}

\def\ppp{\hookrightarrow}
\def\spp{{\overset{s}\hookrightarrow}}
\def\bpp{{\overset{b}\hookrightarrow}}
\def\npp{{\not\hookrightarrow}}
\def\nspp{{\overset{s}{\not\hookrightarrow}}}

\def\wspp{{\overset{w. s}{\hookrightarrow}}}

\newcommand{\pe}{\preccurlyeq}
\newcommand{\npe}{\not\preccurlyeq}
\newcommand{\se}{\succcurlyeq}

\begin{document}
\title[On Stable embeddability of partitions]{On stable embeddability of partitions}
\author{Dongseok KIM}
\address{Department of Mathematics \\ Kyungpook National University \\ Taegu 702-201 Korea }
\email{dskim@math.ucdavis.edu, dongseok86@yahoo.co.kr}
\thanks{}
\author{Jaeun Lee}
\address{Department of Mathematics\\ Yeungnam University\\ Kyongsan, 712-749, Korea }
\email{julee@yu.ac.kr}
\thanks{The first author was supported in part by KRF Grant
M02-2004-000-20044-0. The second author was supported in part by
Com$^2$Mac-KOSEF} \subjclass[2000]{Primary  05A17; Secondary 94A99}
\begin{abstract}
Several natural partial orders on integral partitions, such as the
embeddability, the stable embeddability, the bulk embeddability and
the supermajorization, raise in the quantum computation, bin-packing
and matrix analysis. We find the implications between these partial
orders. For integral partitions whose entries are all powers of a
fixed number $p$, we show that the embeddability is completely
determined by the supermajorization order and we find an algorithm
to determine the stable embeddability.
\end{abstract}

\maketitle

\section{introduction}

A \emph{partition} $\lambda$ is a finite sequence of nonincreasing
positive real numbers, denoted by $\lambda=[\lambda_1, \lambda_2,
\ldots, \lambda_n]$ where $\lambda_i\ge \lambda_j$ for all $i\le j$.
$\lambda_i$ is called an \emph{entry} of $\lambda$. A partition
$\lambda$ is an \emph{integral} partition if all $\lambda_i \in
\mathbb{N}$. Throughout the article, we assume all partitions are
integral unless we state differently. Let $\lambda=[\lambda_1,
\lambda_2, \ldots, \lambda_m]$, $\mu=[\mu_1,\mu_2, \ldots, \mu_n]$
be two partitions. We can naturally define an \emph{addition} of two
partitions, $\lambda + \mu$ by a reordered juxtaposition, a
\emph{product} of two partitions, $\lambda\times\mu$ by
$[\lambda_i\cdot \mu_j]$ and a scalar multiplication,
$\alpha\lambda$ by $[\alpha\cdot \lambda_i]$. We denote
$\overset{n}{\overbrace{\lambda\times\lambda\times
\ldots\times\lambda}}$ by $\lambda^{\times n}$. We recall
definitions of partial orders on partitions. For more terms and
notations, we refer to \cite{bhatia:gtm, Stanley:enumerative2}. A
partition $\lambda$ \emph{supermajorizes} a partition $\mu$, or
$\lambda \se_S \mu$, if for every $x \in \mathbb{N}$

$$\sum_{\lambda_i\ge x} \lambda_i \ge \sum_{\mu_j\ge x} \mu_j.$$
A partition $\lambda$ \emph{embeds} into $\mu$ if there exists a map
$\varphi : \{1, 2, \ldots, m\} \rightarrow \{1, 2, \ldots, n\}$ such
that

$$\sum_{i\in\varphi^{-1}(j)} \lambda_i \le \mu_j$$
for all $j$, denoted by $\lambda\ppp\mu$. This embedding problem can
be interpolated as a \emph{bin-packing problem} by replacing the
entries of a partition $\lambda$ by the sizes of the blocks and the
entries of a partition $\mu$ by the sizes of the bins. It is well
known that the question of whether $\lambda$ embeds into $\mu$ is
computable but NP-hard.

Kuperberg found an interesting embeddability, $\lambda$
\emph{bulk-embeds} into $\mu$, or $\lambda \bpp \mu$, if for every
rational $\epsilon > 0$, there exists an $N$ such that
$\lambda^{\times N} \ppp \mu^{\times
N(1+\epsilon)}$~\cite{Kuperberg:hybrid}. He showed the following
theorem.

\begin{thm}\rm{~\cite{Kuperberg:hybrid}} Let $\lambda$ and $\mu$ are two
partitions, then $\lambda \bpp \mu$ if and only if
$$||\lambda||_p \le ||\mu||_p$$
for all $p \in [1,\infty]$. \label{th:embed}
\end{thm}

He also showed the following implications,

\begin{align}
\lambda \ppp \mu \Longrightarrow\lambda \pe_S \mu
\Longrightarrow \lambda \bpp \mu, \nonumber\\
\lambda\bpp \mu \not\Longrightarrow\lambda \pe_S \mu
\not\Longrightarrow \lambda\ppp\mu. \label{kup}
\end{align}

One can consider a partition as the capacity of a quantum memory
\cite{NV:majorization}. Kuperberg introduced a stable embeddability
in the presence of an auxiliary memory \cite{Kuperberg:hybrid}. A
partition $\lambda$ \emph{stably embeds} into a partition $\mu$ if
there exist a partition $\nu$ such that
$\lambda\times\nu\ppp\mu\times\nu$, denoted by $\lambda\spp \mu$.
Then he asked the relation between the stable embeddability and the
supermajorization order. We answer the question and compare these
embeddabilities in section \ref{compa}. A complete classification of
the stable embeddability is unknown. Since the sizes of the
classical memories are all powers of $2$, it is natural to study the
case all entries of partitions are powers of a fixed positive
integer $p$. For these partitions, we find that the embeddability is
completely determined by the supermajorization order. Also we find
an algorithm to determine the stable embeddability in
section~\ref{stable}.

\section{Comparison of embeddabilities}

For partitions $\lambda, \mu$, we find the following diagram about
the implications of these embeddabilities.

$$
\begin{matrix}
\lambda\ppp\mu & \Longrightarrow & \lambda\pe_S\mu \\
\Downarrow & & \Downarrow \\ \lambda\spp\mu & \Longrightarrow &
\lambda\bpp\mu
\end{matrix}
$$
The converse of all implications are false. We provide
counterexamples in Example~\ref{counterexample}. Moreover, there is
no relation between the stable embeddability and the
supermajorization order, which address the question arose
in~\cite{Kuperberg:hybrid}. For these counterexamples, we need to
show a few facts about these embeddabilities. One can see that if
$\lambda\ppp\mu$, then $||\lambda||_p\le ||\mu||_p$ for all $p\in
[1,\infty].$

\begin{thm} Let $\lambda, \mu$ be partitions. If $\lambda\neq\mu$
and $\lambda \ppp\mu$, then  $||\lambda||_p < ||\mu||_p$ for all
$1<p<\infty$. \label{st}
\end{thm}
\begin{proof}
Let
$$\lambda=[a_1, a_2, \ldots, a_l],~~ \mu
= [b_1, b_2, \ldots, b_m].$$

We will prove it by a contradiction. Suppose $\lambda\ppp \mu$,
$\lambda\neq\mu$ and $||\lambda||_p = ||\mu||_p$ for some
$1<p<\infty$. Then there exists a map $\varphi : \{1, 2, \ldots, l\}
\rightarrow \{1, 2, \ldots, m\}$ presenting the embedding. We divide
cases by the sizes of $l, m$. If $l>m$, then there exist $i_1, i_2$
and $j$ such that $\{i_1, i_2\}\subset \varphi^{-1}(j)$. Since
$$\alpha^p+\beta^p < (\alpha+\beta)^p$$ for all $p> 1$ and nonzero
$\alpha, \beta$, we have

\begin{align}
a_{i_1} +a_{i_2} \le b_j \Longrightarrow a_{i_1}^p+a_{i_2}^p < b_j^p
\Longrightarrow ||\lambda||_p^p < ||\mu||_p^p . \label{inject}
\end{align}

If $l<m$, then there is $k$ such that $\varphi^{-1}(k)=\emptyset$
and hence

\begin{align}
\sum_i (a_{i})^p \le (\sum_j (b_j)^p) - (b_k)^p \Longrightarrow
||\lambda||_p^p < ||\mu||_p^p .\label{surject}
\end{align}

If $l=m$, then there exists $j$ such that $a_k=b_k$ for all $k<j$
and $a_j\neq b_j$ because $\mu^{\times n}\neq\mu^{\times n}$.
Obviously we know $a_j < b_j$. Since $\lambda^{\times n}\ppp
\mu^{\times n},$ either two or more boxes embed into the box of the
size $b_j$ or a part of the box of size $b_j$ has not been used. If
two or more boxes of $\lambda^{\times n}$ embed into the box of the
size $b_j$, then we find a contradiction by equation~\ref{inject}.
If a part of the box of size $b_j$ has not been used, then we find a
contradiction by equation~\ref{surject}. \end{proof}

The following corollary shows the essentiality of $\epsilon$ in
Theorem~\ref{th:embed}.

\begin{cor} Let $\lambda, \mu$ be partitions.
If  $||\lambda||_p = ||\mu||_p$ for some $1<p<\infty$ and
$\lambda\neq \mu$, then $\lambda^{\times n}\npp \mu^{\times n}$ for
all $n$. \label{strict}
\end{cor}

\begin{proof}
Suppose $\lambda^{\times n}\ppp \mu^{\times n}$ for some $n$. Since
$\lambda\neq \mu$, we find $\lambda^{\times n}\neq \mu^{\times n}$.
By Theorem \ref{st} if $\lambda^{\times n}\neq\mu^{\times n}$ and
$\lambda^{\times n}\ppp \mu^{\times n}$ for some $n$, then
$||\lambda^{\times n}||_p < ||\mu^{\times n}||_p$ for all
$1<p<\infty$. But one can observe that for any partition $\lambda$,
$$||\lambda^{\times n}||_p=(||\lambda||_p)^n.$$
Thus we find a contradiction that for all $1<p<\infty$,
$$||\lambda||_p < ||\mu||_p.$$
\end{proof}

\begin{cor} Let $\lambda, \mu$ be two partitions.
If $||\lambda||_p = ||\mu||_p$ for some $1<p<\infty$ and
$\lambda\neq \mu$, then $\lambda\nspp \mu$. \label{cor:proper}
\end{cor}
\begin{proof}
Suppose $\lambda\spp \mu$. There exists a partition $\nu$ such that
$\lambda\times \nu \ppp \mu\times\nu$. Since $\lambda\times \nu \neq
\mu\times\nu$ by Theorem~\ref{st}, for all $1<p<\infty$
$$||\lambda\times \nu||_p  < ||\mu\times\nu||_p.$$
One can easily see that $||\lambda||_p = ||\mu||_p$ for some
$1<p<\infty$ implies that for the same $p$,
$$||\lambda\times \nu||_p = ||\lambda||_p||\nu||_p = ||\mu||_p||\nu||_p =||\mu\times\nu||_p.$$
Therefore, $\lambda\nspp\mu$.
\end{proof}

\begin{exa}
Let $\lambda_1=[2,2,2,2]$, $\lambda_2=[8, 8, 8, 8, 4, 4, 4, 4]$,
$\lambda_3=[4,2,2]$, $\mu_1=[4, \overset{8}{\overbrace{1, 1, \ldots,
1}}]$, $\mu_2=[3,3,3]$, $\mu_3=[16,\overset{16}{\overbrace{2, 2,
\ldots, 2}},\overset{16}{\overbrace{1, 1, \ldots, 1}}]$ and
$\mu_4=[5,3]$. Then $\lambda_1\spp\mu_1$ but $\lambda_1\npe_S\mu_1$
and $\lambda_1\npp\mu_1$. $\lambda_1$$\pe_S$$\mu_2$ but $\lambda_1
\npp \mu_2$. $\lambda_2\bpp\mu_3$ but $\lambda_2\nspp\mu_3$ and
$\lambda_2\npe_S\mu_3$. $\lambda_3\pe_S\mu_4$ but
$\lambda_3\nspp\mu_4$.\label{counterexample}
\end{exa}
\begin{proof}
If we set $\nu=[2,1,1]$, we get
$$\lambda_1\times\nu=[4,4,4,4,\overset{8}{\overbrace{2, 2, \ldots,
2}}]~~ \mathrm{and} ~~
\mu_1\times\nu=[8,4,4,\overset{8}{\overbrace{2, 2, \ldots, 2}},
\overset{16}{\overbrace{1, 1, \ldots, 1}}].$$ Then one can see that
$\lambda_1\spp\mu_1$. Since $$\sum_{(\lambda_1)_i\ge 2}
(\lambda_1)_i =8
> 4 =\sum_{(\mu_1)_j\ge 2} (\mu_1)_j,$$ we see $\lambda_1\npe_S\mu_1$.
It is clear that $\lambda_1\npp\mu_1$. To show $\lambda_2\bpp\mu_3$,
one can check that
$$||\lambda_2||_p \le ||\mu_3||_p$$ for all $p\in[1,\infty]$ and
the equality holds at
$$p=\frac{\mathrm{Ln}(1 + \sqrt{5})}{\mathrm{Ln}(2)}>1.$$
Since $\lambda_2 \neq \mu_3$, we find that $\lambda_2 \nspp \mu_3$
by Corollary~\ref{cor:proper}. Clearly $\lambda_3\pe_S\mu_4$.
Suppose $\lambda_3\spp\mu_4$, there exists a partition
$\nu=[\nu_1,$$ \nu_2,$$ \ldots,
$$\nu_n]$ such that $\lambda_3\times\nu \ppp \mu_4\times\nu$. Let
$p$ be the power of $2$ in the prime factorization of the greatest
common divisor $(\nu_1, \nu_2, \ldots, \nu_n)$ of $\nu_1, \nu_2,
\ldots, \nu_n$. First we looks at entries of $\lambda_3\times\nu$,
all these entries are multiples of $2^{p+1}$. Since
$||\lambda_3\times\nu||_1=||\mu_4\times\nu||_1$, there will be no
space in $\mu_4\times \nu$ which was not used in the embedding, $i.
e.,$ for all $j$,
$$\sum_{i\in\varphi^{-1}(j)} (\lambda_3\times\nu)_i = (\mu_4\times\nu)_j.$$
Therefore, all entries of $\mu_4\times\nu$ have to be multiples of
$2^{p+1}$. Since all entries of $\mu_4$ are odd numbers, $\nu_i$ has
to be a multiple of $2^{p+1}$ and so does the greatest common
divisor of $\nu_1, \nu_2, \ldots, \nu_n$. It contradicts the
hypothesis of $p$. All others should be straightforward.
\end{proof}
\label{compa}

\section{Stable embeddability}

Let $\lambda, \mu$ be two partitions. Let us consider the following
algorithm which is called a \emph{ first fit}
algorithm~\cite{AM:firstrandom}. From $\lambda_1$ of $\lambda$,
place it to any entry of $\mu$ in which it fits. Then repeat this
step for $\lambda_2$ and so on. Usually this is not an efficient
algorithm~\cite{johnson:bestvsfirst}. It is obvious that if the
first fit algorithm works, then $\lambda\ppp \mu$. The converse is
not true in general. But with some conditions on $\lambda$ we can
show that it determines the embeddability of $\lambda$ into $\mu$.

\begin{thm}
Let $\lambda=[\lambda_1,\lambda_2, \ldots, \lambda_s], \mu =[\mu_1$,
$\mu_2$, $\ldots $, $\mu_t]$ be partitions with
$\lambda_i|\lambda_j$ for all $i\ge j$. If $\lambda\ppp \mu$, then
the first fit algorithm works. \label{dividable}
\end{thm}
\begin{proof}
Let us induct on $s$. It is trivial for $s=1$ because this is the
first step of the algorithm. Suppose this is true for $s=n$, we look
at the case $\lambda=[\lambda_1,\lambda_2, \ldots,\lambda_{n+1}]$.
Since $\lambda\ppp \mu$, there is a map $\varphi : \{ 1, 2, \ldots ,
n+1\} \rightarrow \{1, 2, \ldots, t\}$ which represents the
embedding of $\lambda$ into $\mu$, let us denote $\varphi(1)=j$.
Then we will construct another embedding representing map $\psi$
after we decide where we put $\lambda_1$, say $\mu_k$, $i.e.,$ $\psi
(1)=k$. To construct $\psi$, let us compare $\varphi(1)$ and
$\psi(1)$. If $\varphi(1)=\psi(1)$, we pick $\psi = \varphi$. If
$\varphi(1)=j\neq k=\psi(1)$, first we need to prove that there
exists a subset $P$ of $\varphi^{-1}(k) =\{\lambda_{i_1},
\lambda_{i_2}, \ldots , \lambda_{i_l}\}$ such that
$$\sum_{j\in P} \lambda_j \le \lambda_1 \hskip 1cm \mathrm{and} \hskip 1cm
\sum_{j\in P^c} \lambda_j \le \mu_k-\lambda_1.$$ First we divide
cases by the sizes of $\lambda_1$ and $\lambda_2$. If
$\lambda_1=\lambda_2$, then we pick $P=\{2\}$. If $\lambda_1
\neq\lambda_2$, then $\lambda_1>\lambda_2$. If
$$\sum_{j=2}^{l} \lambda_j \le \lambda_1,$$ we can pick $P=\{2, 3,
\ldots , l\}$. Otherwise, there exist an integer $m$ such that
$$ \sum_{j=2}^{m} \lambda_j\le \lambda_1 < \sum_{j=2}^{m+1}\lambda_j
= (\sum_{j=2}^{m}\lambda_j) + \lambda_{m+1}.$$ If we divide by
$\lambda_{m+1}$, we have
$$ \sum_{j=2}^{m} \frac{\lambda_j}{\lambda_{m+1}}\le
\frac{\lambda_1}{\lambda_{m+1}} <
\sum_{j=2}^{m}\frac{\lambda_j}{\lambda_{m+1}} + 1.$$ Since
$\lambda_i|\lambda_j$ for all $i\ge j$, all these three numbers are
integers, so the first two have to be the same. We choose $P=\{2, 3,
\ldots , m\}$. Once we have such a $P$, we can define
$$\psi(i) = \left\{
\begin{array}{cl}
 k  & ~~\mathrm{if}~~  i=1, \\
 j  & ~~\mathrm{if}~~ i \in \varphi^{-1}(k)\cap P, \\
  \varphi (i) & ~~\mathrm{if}~~ i\not\in\varphi^{-1}(k)
      \cup\{1\} ~~\mathrm{or}~~ i\in \varphi^{-1}(k)\cap P^{c}.
\end{array} \right.
$$
Then $\psi|_{\tilde\lambda}$ shows $\tilde\lambda=[\lambda_2,\ldots
\lambda_{n+1}] \ppp \tilde\mu=[\mu_1, \ldots,
\mu_k-\lambda_1,\ldots, \mu_t]$. By the induction hypothesis, the
first fit algorithm works.
\end{proof}

Let $\mathcal{P}$ be the set of all partitions whose entries are all
powers of a fixed number $p$.  For these partitions, we can show
that the supermajorization completely determine the embeddability.
Instead of the standard notation, we can use
$$\lambda=[a_0, a_1, a_2, \ldots, a_s]_p$$
where $a_i$ is the number of entries $p^i$.

\begin{thm} Let $\lambda, \mu$ be partitions in $\mathcal{P}$.
$\lambda\ppp\mu$ if and only if
$\lambda\pe_S \mu$. \label{superpp}
\end{thm}
\begin{proof}
We only need to show that if $\lambda\pe_S\mu$, $\lambda\ppp\mu$
because of equation \ref{kup}. Suppose $\lambda\pe_S\mu$. Let
$$\lambda=[ a_0, a_1,a_2,
\ldots, a_s]_p, ~~ \mu=[b_0, b_1, b_2, \ldots, b_t]_p.$$ Without
loss of generality, we assume $a_s\neq 0\neq b_t$. Obviously $s\le
t$. We induct on the number of the boxes of $\lambda$, say $k$. If
$k=1$, then $a_s=1$
$$
1p^s=\lambda_{\ge p^s} \le \mu_{\ge p^s}=\sum_{j=s}^{t}b_jp^{j}
$$
implies $\lambda\ppp\mu$. For nonzero $a_s$, we pick a box of size
$p^s$, put it into a box of size $p^t$ in $\mu$. Then for $\lambda$
we subtract $1$ from $a_s$ and for $\mu$, we subtract $1$ from $b_t$
and distribute the reminder of $p^t-p^s$  in base $p$ into $\mu$.
One can observe that all these numbers which have been distributed
are bigger than or equal to $p^s$. Thus resulting partitions still
have the same supermajorization order. By the induction hypothesis,
we find an embedding of $\lambda'=[a_0, a_1, \ldots, a_{s-1},
a_s-1]_p$ into $\mu'=[b_0', b_1', \ldots, b_{t-1}', b_t-1]_p$. But,
it is easy to recover an embedding of $\lambda$ into $\mu$.
\end{proof}

Now we look the stable embeddability for partitions in
$\mathcal{P}$.

\begin{thm}
Let $\lambda, \mu$ be partitions in $\mathcal{P}$. If $\lambda\spp
\mu$, then there exists a partition $\nu$ in $\mathcal{P}$ such that
$\lambda\times\nu \ppp \mu\times\nu$. \label{Cpower}
\end{thm}

\begin{proof}
Suppose $\lambda\spp \mu$, then there is a partition $\nu$ such that
$\lambda\times\nu \ppp \mu\times\nu$ and $\nu=[c_1, c_2, \ldots,
c_k].$ We can uniquely rewrite $c_j$ in the base $p$ such as
$$c_j= c_{j,0} p^0 + c_{j,1}p^1 + c_{j,2} p^{2} + \ldots + c_{j,l(j)}
p^{l(j)}$$ where $c_{j,i}$ are nonnegative integers less than $p$
and $ c_{j,l(j)}\neq 0$. Using these expressions we can subdivide
$\nu$ to get a refinement $$\tilde\nu=[\sum_{j} c_{j,0}, \sum_{j}
c_{j,1}, \ldots , \sum_j c_{j,i}, \ldots, \sum_{j} c_{j, t}]_p$$
where the sum runs over all nonzero $c_{j,i}$ for each $i$. If the
boxes $\sum[ c_{i_k}\times p^{j_k}]$ of $\lambda\otimes\nu$ were
embedded into $[c_m\times p^{m'}]$ in $\mu\otimes\nu$, We can show
that the refinement of $\sum[ c_{i_k}\times p^{j_k}]$ can be
embedded in the refinement of  $[c_m\times p^{m'}]$. Precisely if
$$p^{j_1}c_{i_1}+ p^{j_2}c_{i_2} + \ldots + p^{j_n}c_{i_n} \le
p^{m'}c_m$$ where $j_1\le j_2 \le \ldots \le j_n$, $c_{j_t}\neq 0$
and
$$c_{i_t}=c_{i_t,0}p^{0} + c_{i_t,1}p^{1} + \ldots + c_{i_t,l(i_t)}
p^{l(i_t)}$$ for all $t$, then
$$\sum_{\alpha =1}^{n}\sum_{\beta =0}^{l(\beta)}
[c_{i_{\alpha},\beta}\times p^{i_{\alpha}+\beta}]\ppp
\sum_{\gamma}^{l(m)}[c_{m,\gamma}\times p^{m'+\gamma}].$$

First we look at the case, $n = 1$. If $p^{j_1}c_{i_1} \le p^{m'}
c_m$, one can easily see that
$$\sum_{\beta}[c_{i_1,\beta}\times p^{j_1+\beta}]\pe_S
\sum_{\gamma}[c_{m,\gamma}\times p^{m'+\gamma}]$$ because we are
comparing two integers in base $p$. By Theorem~\ref{superpp},
$$\sum_{\beta}[c_{i_1,\beta}\times p^{j_1+\beta}]\ppp
\sum_{\gamma}[c_{m,\gamma}\times p^{m'+\gamma}].$$

For the case $n>1$, we look at the integer $$\sum_{\alpha
=1}^{n}\sum_{\beta =0}^{l(\beta)} c_{i_{\alpha},\beta}\times
p^{i_{\alpha}+\beta}$$ as a sum of integers $$\sum_{\beta
=0}^{l(\beta)} c_{i_{\alpha},\beta}\times p^{i_{\alpha}+\beta}$$ in
base $p$. Then this returns to the case $n=1$. If we keep on
tracking the addition, we can recover the embedding of
$$\sum_{\alpha =1}^{n}\sum_{\beta =0}^{l(\beta)}
[c_{i_{\alpha},\beta}\times p^{i_{\alpha}+\beta}]\ppp
\sum_{\gamma}^{l(m)}[c_{m,\gamma}\times p^{m'+\gamma}].$$

Moreover, this process does not involve with other terms. Therefore,
we can rewrite $\nu$ as the shape we desired.
\end{proof}

\begin{cor}
Let $\lambda=[a_i]_p, \mu=[b_i]_p$ be partitions in $\mathcal{P}$.
If $0 \le a_i\le b_i (0\le b_i < a_i)$, we have two new partitions
$\tilde\lambda,$$ \tilde\mu$ which are obtained from $\lambda,$
$\mu$ by replacing $a_i, b_i$ by $0$ and from $\mu, \lambda$ by
replacing by $b_i-\mathrm{Min}\{ a_i, b_i\}(a_i-\mathrm{Min}\{ a_i,
b_i\},$ $respectively)$. Then $\lambda\spp \mu$ if and only if
$\tilde\lambda\spp \tilde\mu$. \label{shape}
\end{cor}

\begin{proof}
We assume $a_i\le b_i$ for a fixed $i$. Suppose $\lambda\spp \mu$.
By Theorem~\ref{Cpower}, we can find

$$\nu=[\nu_0, \nu_{1}, \ldots, \nu_n]_p$$
such that all entries of $\nu$ are all powers of a fixed number $p$
and $c_k$ is the number of the boxes of size $p^k$. Now
$\lambda\otimes\nu\ppp \mu\otimes\nu$ and $\lambda\otimes\nu,
\mu\otimes\nu$ satisfy the hypothesis of Theorem~\ref{dividable}, we
can use the first fit algorithm. We put all boxes whose sizes are
bigger than $p^i\times p^n$. Then we consider $a_i\cdot c_n$ boxes
of size $p^i\times p^n$ in $\lambda\otimes\nu$. But none of boxes of
size $p^i\times p^n$ in $\mu\otimes\nu$ were used in the previous
steps, we can put these into boxes of size $p^i\times p^n$ of
$\mu\otimes\nu$. Then we finish the rest of boxes of size $p^{n+i}$.
Then we repeat the same process to the next $a_i\cdot c_{n-1}$ boxes
of size $p^i\times p^{n-1}$ in $\lambda\otimes\nu$. This embedding
keeps all boxes $p^{i}\otimes \nu$ of $\lambda\otimes\nu$ into
$p^{i}\otimes \nu$ of $\mu\otimes\nu$. Thus, $\lambda\otimes\nu\ppp
\mu\otimes\nu$. The converse is obvious.
\end{proof}

\subsection{An algorithm to determine the stable embeddability}
We introduce an algorithm to decide $\nu$. Let $\lambda, \mu$ be
partitions in $\mathcal{P}$, $i.e.,$ $\lambda=[a_0, a_1, \ldots,
a_n]_p$ and $\mu=[b_0, b_{1}, \ldots, b_m]_p,$ where $a_i, b_i$ are
the number of boxes of size $p^i$ in $\lambda, \mu$ respectively. By
Theorem~\ref{superpp}, we can decide whether $\lambda$ can be
embedded in $\mu$ or not. Before we apply the algorithm, we modify
the shape of $\lambda, \mu$ by Corollary~\ref{shape} such that none
of $a_i, b_i$ are nonzero simultaneously. If $a_n\neq 0\neq b_m$ and
$m<n$, $\lambda$ can not be stably embedded into $\mu$. For
convenience, we will assume $p$ is $2$, $\nu$ is a rational
partition whose entries are of non-positive powers of $2$ and $c_k$
is the number of boxes in $\nu$ of the size $2^{-k}$.

Initially, we will start $c_0=1$. There are $a_n$ boxes of size
$2^n$ in $\lambda\times[c_0 \times 1]$ and none of blocks of size
$2^n$ in $\mu\times[c_0\times 1]$. But there are rooms for

$$b_m \times 2^{m-n} +b_{m-1} \times 2^{m-n-1} + \ldots + b_{n+1}\times 2$$
many boxes of size $2^n$ in $\mu\times[c_0\times 1]$. If

$$b_m \times
2^{m-n} + b_{m-1} \times 2^{m-n-1} + \ldots + b_{n+1}\times 2 \ge
a_n,$$ we set $c_1$ to zero and keep the difference for the next
step, say $M$. Otherwise we set $$c_1 =\lceil \frac{a_n - (b_m
\times 2^{m-n} + b_{m-1} \times 2^{m-n-1} + \ldots + b_{n+1}\times
2)}{b_m} \rceil$$ and $M=0$, where $\lceil x \rceil$ is the smallest
natural number which is bigger than or equal to $x$. Then we look at
$\lambda\times[c_0\times 1, c_1\times\frac{1}{2}],
\mu\times[c_0\times 1,c_1\times\frac{1}{2}]$. We have $a_{n-1}\cdot
c_0 + a_n\cdot c_1$ many boxes of size $2^{n-1}$ in
$\lambda\times[c_0\times 1]$, then we compare it with

$$
2\times M+ c_1\times ( b_m\times 2^{m-n} + b_{m-1} \times 2^{m-n-2}
+ \ldots + b_{n+1}\times 2) + b_{n-1}\cdot c_0
$$
and we repeat exactly the same process. For $N\ge m$, we find $c_N$
by comparing two terms $$\alpha= b_{m-1}\times c_{N+m-1} +
b_{m-2}\times c_{N+m-2} + \ldots + b_0\times c_{N} + 2\times M$$ and
$$\beta=a_{n}\times c_{N+n} + a_{n-1}\times c_{N+n-1} + \ldots +
a_{0}\times c_{N}$$ because these numbers count exactly how many
blocks of size $2^{-N}$ in the product
$$\lambda\times[c_0\times 1, c_1\times\frac{1}{2},
\ldots, c_N\times\frac{1}{2^{N}}]$$ and $$\mu\times[c_0\times1,
c_1\times\frac{1}{2}, \ldots, c_N\times\frac{1}{2^N}]$$ where $M$ is
the number of boxes that were left in the previous step. Then
$C_{N+1}$ is $\lceil (\beta-\alpha)/b_m \rceil$ if $\beta-\alpha > 0
($and set $M=0)$ $0$ otherwise $($set $M= \alpha-\beta,$
$respectively)$. Then we compare the next biggest boxes. We stop if
we get $n$ consecutive $0$'s for $c_i$. Let $N$ be the largest
integer that $c_N$ is non-zero. We repeat the process starting
$c_0=(b_m)^{N+1}$. One can easily see that we no longer have to use
$\lceil \rceil$ because $(b_m)^{N+1-k}|c_k$ for all $0\le k \le
N+1$. Finally we multiply $2^N$ to make $\nu$ an integral partition.

To compare the optimality of such $\nu$'s, we define the
\emph{length} of $\nu=[c_0, c_1, \ldots, c_n]_p$ to be $n+1$ where
$c_0\neq 0 \neq c_n$. From the given $\lambda, \mu$ we collect all
possible $\nu \in \mathcal{P}$ and $\lambda\otimes \nu\ppp
\mu\otimes \nu$, say $\tt$. Then we define a partial order on $\tt$
by a lexicographic order,
$$( \mathrm{length}\hskip .2cm \mathrm{of}\hskip .2cm \widehat\lambda(\nu),
\frac{c_1}{c_0}, \ldots, \frac{c_n}{c_0}).$$ Moreover, $\tt$ is
closed under an addition, a tensor and a scalar multiplication.

\begin{thm}
1) The algorithm stops at finite time if and only if
$\lambda\spp\mu$.

2) Let $D$ be a partition which is obtained from the algorithm. Then
$D$ is a minimal element with respect to the partial order we
defined on $\tt$. \label{optimal}
\end{thm}

\begin{proof}
We want to show that if $\lambda\spp\mu$, then the algorithm must
stop at finite steps and the one we find by the algorithm has the
smallest length. Since $\lambda\spp\mu$, $\tt$ is nonempty and we
find a minimal element in $\tt$, say $$\ttt=[t_0, t_1, \ldots,
t_l]_p.$$ First we assume the existence, $i. e.,$ the algorithm
gives us an integral partition
$$\nu=[c_0, c_1, \ldots, c_m]_p.$$ By the minimality,
we have $l\le m$. But the process itself provides us $l\ge m$. We
compare $$c_0\ttt=[c_0\cdot t_0, c_0\cdot t_1, \ldots, c_0\cdot
t_l]_p$$ and $$t_0 \nu=[t_0\cdot c_0, t_0\cdot c_1, \ldots, t_0\cdot
c_m]_p.$$ Suppose $c_o\ttt\neq t_o\nu$. There is a $j$ such that
$$c_0\cdot t_j < t_0\cdot c_j.$$ But this obviously contradicts the
process of the algorithm. Therefore, $c_0\ttt = t_0\nu$ and it does
also prove the existence.
\end{proof}

\label{stable}

\section{Discussions}
\label{discussion}

\subsection{Algebraic embeddabilities.}
Let $\aa$ be a finite dimensional semisimple algebra over an
algebraically closed field $K$. By a simple application of
Webberburn-Artin theorem, we can decompose $\aa$ into a direct sum
of matrix algebras. From a direct sum of matrix algebras $\aa$, we
can find a unique integral partition $\lambda$, denoted by
$\lambda(\aa)$. For an integral partition $\lambda$, one can assign
a direct sum of matrix algebras
$$\aa(\lambda) = \bigoplus_{i=1}^{m} {\mathcal M}_{\lambda_i},$$
where ${\mathcal M}_{\lambda_i}$ is the set of all $\lambda_i$ by
$\lambda_i$ matrices over $K$. For integral partitions, one can see
that $\lambda\ppp\mu$ if and only if $\aa(\lambda)$ embeds into
$\aa(\mu)$ as $K$ algebras. All other partial orders can be
naturally defined for a direct sum of matrix algebras. The question
of the embeddability between algebraic objects such as groups,
rings, modules and etc, is a long standing difficult question. For
some algebraic objects such as sets, vector spaces, the question is
straightforward. The embeddability between the modules over a
complex simple Lie algebra is completely determined by
Littlewood-Richardson formula and Schur's lemma. Authors have made a
few progress on stable embeddability, the product is replaced by the
tensor product, between the modules over a complex simple Lie
algebra \cite{lie}. The stable embeddability between other algebraic
objects should be an interesting question.

\subsection{Analytic embeddabilities.}

Let $\lambda, \mu$ be partition in $\mathcal{P}$. The algorithm we
defined in section \ref{stable} brings us a new embeddability,
$\lambda$ \emph{weakly stably embeds} into $\mu$, denoted by
$\lambda\wspp\mu$, if there exists a rational partition $\nu$ of
infinite length such that all entries of $\nu$ are nonpositive
powers of the fixed number $p$ and
$$\sum_{i=0}^{\infty} c_i p^{-i} <\infty,$$
where $c_i$ is the number of the entries $p^{-i}$. One can see that

\begin{eqnarray}
\begin{matrix}
\lambda\spp\mu & \Longrightarrow & \lambda\wspp\mu  & \Longrightarrow & \lambda\bpp\mu \\
& &  & & \Updownarrow \\
& & ||\lambda||_p < ||\mu||_p, ~ \forall p\in(1,\infty) &
\Longrightarrow & ||\lambda||_p \le ||\mu||_p,~ \forall
p\in[1,\infty]
\end{matrix}
\label{ws}
\end{eqnarray}

It is not known that the converses of the first row of equation
\ref{ws} are true or not for partitions in $\mathcal{P}$. Authors
have written a program that performs the algorithm described in
section~\ref{stable} to see $||\lambda||_p < ||\mu||_p, ~ \forall
p\in(1,\infty)$ and equality holds for $p=1$ and $\infty$ implies
$\lambda\spp\mu$. We have not found any answer yet.

\bibliographystyle{amsplain}

\end{document}